\documentclass[11pt]{article}

\usepackage{amssymb}
\usepackage{amsfonts}
\usepackage{eucal}
\usepackage{soul}
\usepackage{amsmath}
\newtheorem{theorem}{Theorem}[section]
\newtheorem{proposition}[theorem]{Proposition}
\newtheorem{definition}[theorem]{Definition}
\newtheorem{corollary}[theorem]{Corollary}

\newtheorem{remark}[theorem]{Remark}

\newenvironment{proof}[1][Proof]{\textbf{#1.} }{\ \rule{0.5em}{0.5em}}
\input pictex.sty

\usepackage{cancel}
\usepackage{xcolor}

\usepackage{amsmath,color,ulem}
\newcommand{\stkout}[1]{\ifmmode\text{\sout{\ensuremath{#1}}}\else\sout{#1}\fi}

\DeclareMathOperator\E{E}
\DeclareMathOperator\Var{Var}
\DeclareMathOperator\Cov{Cov}

\begin{document}
	\title{Doob-type optional sampling theorems for demimartingales with applications to associated sequences}
	\author{Milto Hadjikyriakou\footnote{School of Sciences, University of Central Lancashire, Cyprus campus, 12-14 University Avenue, Pyla, 7080 Larnaka, Cyprus (email:mhadjikyriakou@uclan.ac.uk).}~~ and B.L.S. Prakasa Rao \footnote{CR RAO Advanced Institute of Mathematics, Statistics and Computer Science, Hyderabad 500046, India (e-mail: blsprao@gmail.com).}}	
	\maketitle	
	
	\begin{abstract}
	We establish optional sampling inequalities for demimartingales and demisubmartingales under suitable monotonicity assumptions on the stopping rule. First, we establish comparison inequalities for processes stopped at bounded stopping times and identify broad classes of threshold-type stopping times for which the required assumptions are naturally satisfied. We then prove Doob-type optional sampling inequalities for possibly unbounded stopping times under standard integrability conditions. As applications, we derive maximal inequalities for demi(sub)martingales, obtain an almost sure bound for the maximal negative excursion of partial sums of positively associated random variables, and establish generalized Wald-type inequalities, including nonlinear convex-transform inequalities for associated sequences.
	
	\medskip
	
	
	\end{abstract}
	
	\section{Introduction}
	
	Demimartingales were originally introduced in \cite{NW1982} as a tool to study dependent structures beyond the scope of classical martingales, particularly in the context of positively associated sequences (see \cite{EPW1967}). The relevant definitions are provided below.
	
	\begin{definition}
		A finite collection of random variables $X_1, \ldots, X_n$ is said to be (positively) associated if
		$$
		\operatorname{Cov}\left(f\left(X_1, \ldots, X_n\right), g\left(X_1, \ldots, X_n\right)\right) \geq 0
		$$
		for any componentwise nondecreasing functions $f, g$ on $\mathbb{R}^n$ such that the covariance is defined. An infinite collection is associated if every finite subcollection is associated.
	\end{definition}
	
	\begin{definition}
		A sequence of $L^1$ random variables $\left\{S_n, n \in \mathbb{N}\right\}$ is called a demimartingale if for all $j=1,2, \ldots$
		$$
		\E\left[\left(S_{j+1}-S_j\right) f\left(S_1, \ldots, S_j\right)\right] \geq 0
		$$
		for all componentwise nondecreasing functions $f$ whenever the expectation is defined. Moreover, if $f$ is assumed to be nonnegative, the sequence $\left\{S_n, n \in \mathbb{N}\right\}$ is called a demisubmartingale.
	\end{definition}
	
	For more on demimartingales, the interested reader may refer to \cite{C2000}, \cite{C2003}, \cite{DSHY2014}, \cite{HPrR2025} and the monograph of Prakasa Rao \cite{PR2012} while for more on associated random variables, we refer to \cite{O2012}. 
	
	\begin{remark}
		\label{rem1}
		Consider a sequence of positively associated random variables $(X_n)_{n\geq 1}$ and let $S_n = \sum_{i=1}^nX_i$. Then, for any componentwise nondecreasing function $f$, we have:
		\[
		\E[(S_{n+1}- S_n)f(S_1, \ldots,S_n)] = \E[X_{n+1}f_1(X_1, \ldots,X_n)] \geq \E X_{n+1} \E[f_1(X_1,\ldots,X_n)]
		\]
		due to the association property. Observe that:
		\begin{enumerate}
			\item [(a)] if $\E X_i = 0$ for all $i$, then $S_n$ forms a demimartingale sequence 
			\item [(b)]  if $\E X_i \geq 0$ for all $i$ and $f$ is also assumed to be nonnegative then $S_n$ is a demisubmartingale.
		\end{enumerate}
		Counterexamples available in the literature  (see for example \cite{MH2010} or \cite{PR2012}) prove that not all demimartingale sequences have associated demimartingale differences, i.e. demimartingales form a class wider than the class of partial sums of zero mean positively associated random variables.
	\end{remark}

	\begin{remark}
		\label{rem2}Martingales satisfy the definition of demimartingales when considered with their standard filtration. Thus, every martingale with the natural choice of $\sigma$-algebras is also a demimartingale. However, the converse does not hold in general (see \cite{MH2010} or \cite{PR2012} for counterexamples).
	\end{remark}

		The examples discussed above illustrate two classical sources of demimartingales, namely partial sums of positively associated random variables and martingales. The class of demimartingales, however, is considerably broader. The following propositions present three families of demimartingales arising from latent-factor models, dynamic observation systems, and multiplicative growth processes. Their practical interpretations will be discussed after each proposition. They further demonstrate that demimartingales need not be martingales, nor do they necessarily possess associated increments.
		
		\begin{proposition}
			\label{prop:gaussian-latent}
			Let
			\[
			S_n = \beta_n Z + \varepsilon_n, \quad n \geq 1,
			\]
			where $Z, \varepsilon_1, \varepsilon_2, \ldots$ are independent centred normally distributed
			random variables with $\Var(Z) = \sigma_1^2 > 0$ and $\Var(\varepsilon_n) = \sigma_2^2 > 0$,
			and let $(\beta_n)_{n \geq 1}$ be a nondecreasing sequence of nonnegative numbers. If
			\[
			(\beta_{j+1} - \beta_j)\beta_j \sigma_1^2 \geq \sigma_2^2 \quad \text{for all } j \geq 1,
			\]
			then $\{S_n,\, n \geq 1\}$ is a demimartingale.
		\end{proposition}
		\begin{proof}
			Fix $j \geq 1$ and set $\Delta_j = \beta_{j+1} - \beta_j$, so that
			\[
			S_{j+1} - S_j = \Delta_j Z + \varepsilon_{j+1} - \varepsilon_j.
			\]
			Each component of the vector $(S_1, \ldots, S_j, S_{j+1} - S_j)$ is a linear combination
			of the independent Gaussian variables $Z, \varepsilon_1, \ldots, \varepsilon_{j+1}$, so the
			vector is jointly Gaussian. For Gaussian vectors, nonnegative pairwise covariances imply
			association, and it therefore suffices to verify these. For $1 \leq i, k \leq j$,
			\[
			\Cov(S_i, S_k) = \beta_i \beta_k \sigma_1^2 + \sigma_2^2 I{\{i=k\}} \geq 0,
			\]
			since $\beta_i, \beta_k \geq 0$. For $1 \leq i < j$,
			\[
			\Cov(S_{j+1} - S_j,\, S_i) = \Delta_j \beta_i \sigma_1^2 \geq 0,
			\]
			where the other terms vanish because $\varepsilon_i$ is independent of both $\varepsilon_j$
			and $\varepsilon_{j+1}$. Finally,
			\[
			\Cov(S_{j+1} - S_j,\, S_j) = \Delta_j \beta_j \sigma_1^2 - \sigma_2^2 \geq 0
			\]
			by assumption. Thus $(S_1, \ldots, S_j, S_{j+1} - S_j)$ is associated, and for every
			componentwise nondecreasing $g : \mathbb{R}^j \to \mathbb{R}$,
			\[
			\E\bigl[(S_{j+1} - S_j)\,g(S_1, \ldots, S_j)\bigr]
			= \Cov\bigl(S_{j+1} - S_j,\, g(S_1, \ldots, S_j)\bigr) \geq 0,
			\]
			where the first equality uses $\E(S_{j+1} - S_j) = 0$.
		\end{proof}
		\begin{remark}
			\label{rem:not-martingale-1}
			The process $\{S_n\}$ is not a martingale with respect to its natural filtration
			$\mathcal{F}_n = \sigma(Z, \varepsilon_1, \ldots, \varepsilon_n)$, since
			\[
			\E(S_{n+1} \mid \mathcal{F}_n) = \beta_{n+1} Z \neq \beta_n Z + \varepsilon_n = S_n.
			\]
		
		\medskip
			
			\noindent The increment sequence $X_n = S_{n+1} - S_n$ need not be associated either. Writing
			$X_n = \Delta_n Z + \varepsilon_{n+1} - \varepsilon_n$ and
			$X_{n+1} = \Delta_{n+1} Z + \varepsilon_{n+2} - \varepsilon_{n+1}$, one computes
			\[
			\Cov(X_n, X_{n+1}) = \Delta_n \Delta_{n+1} \sigma_1^2 - \sigma_2^2.
			\]
			If $\Delta_n \Delta_{n+1} \sigma_1^2 < \sigma_2^2$ for some $n \geq 1$, this covariance is
			negative. Since $(X_n, X_{n+1})$ is Gaussian, association would require all pairwise
			covariances to be nonnegative, so the increment sequence is not associated (see \cite{P1982}).
		\end{remark}
		
		\medskip
		
		\begin{remark}
			\label{rem:interpretation-1}
			The process $S_n = \beta_n Z + \varepsilon_n$ models a sequence of noisy observations of
			an unobservable latent state $Z$, where $\varepsilon_n$ represents measurement error and
			$\beta_n$ describes the sensitivity of the observation mechanism. The condition
			$\beta_{n+1} \geq \beta_n$ corresponds to progressively more informative observations over
			time. A natural stopping time is
			\[
			\tau_b = \inf\{n \geq 1 : S_n > b\},
			\]
			the first time the observed signal exceeds a prescribed threshold. Such models arise naturally in signal detection, quality control, monitoring systems, and repeated measurements of latent characteristics. In these settings, $S_{\tau_b}$
			represents the observed signal level at the time an alarm or intervention is triggered. Optional sampling results are particularly useful since they provide information about the stopped process without requiring an explicit analysis of the distribution of $\tau_b$.
		\end{remark}
		
		\medskip
		
		\noindent The previous proposition corresponds to a static latent-factor model. The next result extends this framework to a dynamic latent state evolving as a Gaussian random walk observed with noise.
		
		\medskip
		
		\begin{proposition}
			\label{prop:random-walk-factor}
			Let $M_n = \sum_{k=1}^n \xi_k$, $n \geq 1$, where $\xi_1, \xi_2, \ldots$ are independent
			centred normal random variables with $\Var(\xi_k) = \sigma_1^2 > 0$. Let
			$\varepsilon_1, \varepsilon_2, \ldots$ be independent centred normal random variables,
			independent of $\{\xi_k\}$, with $\Var(\varepsilon_n) = \sigma_2^2 > 0$. Define
			\[
			S_n = \beta_n M_n + \varepsilon_n, \quad n \geq 1,
			\]
			where $\{\beta_n\}$ is a nondecreasing sequence of nonnegative numbers satisfying
			\[
			(\beta_{n+1} - \beta_n)\beta_n n\,\sigma_1^2 \geq \sigma_2^2, \quad n \geq 1.
			\]
			Then $\{S_n,\, n \geq 1\}$ is a demimartingale.
		\end{proposition}
		\begin{proof}
			Fix $j \geq 1$. Since $M_{j+1} = M_j + \xi_{j+1}$,
			\[
			S_{j+1} - S_j = \Delta_j M_j + \beta_{j+1}\xi_{j+1} + \varepsilon_{j+1} - \varepsilon_j.
			\]
			The vector $(S_1, \ldots, S_j, S_{j+1} - S_j)$ is jointly Gaussian, so again it suffices
			to check pairwise covariances. For $1 \leq i, k \leq j$,
			\[
			\Cov(S_i, S_k) = \beta_i \beta_k \min(i,k)\,\sigma_1^2 + \sigma_2^2 I{\{i=k\}} \geq 0.
			\]
			For $1 \leq i < j$,
			\[
			\Cov(S_{j+1} - S_j,\, S_i) = \Delta_j \beta_i\, i\,\sigma_1^2 \geq 0,
			\]
			and for $i = j$,
			\[
			\Cov(S_{j+1} - S_j,\, S_j) = \Delta_j \beta_j\, j\,\sigma_1^2 - \sigma_2^2 \geq 0
			\]
			by assumption. Hence $(S_1, \ldots, S_j, S_{j+1} - S_j)$ is associated, and the
			demimartingale property follows as in Proposition~\ref{prop:gaussian-latent}.
		\end{proof}
		
		\medskip

		\begin{remark}
			\label{rem:not-martingale-2}
			Consider the natural filtration $\mathcal{F}_n = \sigma(\xi_1,\ldots,\xi_n, \varepsilon_1, \ldots, \varepsilon_n)$. Then, $\E(S_{n+1} \mid \mathcal{F}_n) = \beta_{n+1} M_n \neq S_n$ i.e. $\{S_n\}$ is not a
			martingale. The increment covariance is
			\[
			\Cov(X_n, X_{n+1}) = \Delta_{n+1}(\Delta_n n + \beta_{n+1})\sigma_1^2 - \sigma_2^2,
			\]
			which is negative whenever $\Delta_{n+1}(\Delta_n n + \beta_{n+1})\sigma_1^2 < \sigma_2^2$.
			Under this condition, $(X_n, X_{n+1})$ has a negative pairwise covariance and, since the
			pair is Gaussian, the increment sequence is not associated (see \cite{P1982}).
		\end{remark}

		\begin{proposition}
			\label{prop:multiplicative}
			Let $R_1, R_2, \ldots$ be independent nonnegative random variables with
			$\E(R_i) = \mu > 1$ for all $i \geq 1$, and assume all expectations below exist. Define
			\[
			P_n = \prod_{i=1}^n R_i, \qquad S_n = P_n - \E(P_n) = P_n - \mu^n.
			\]
			Then $\{S_n,\, n \geq 1\}$ is a demimartingale.
		\end{proposition}
		\begin{proof}
			Fix $j \geq 1$ and let $g : \mathbb{R}^j \to \mathbb{R}$ be componentwise nondecreasing
			with all expectations finite. Since $P_{j+1} = P_j R_{j+1}$,
			\[
			S_{j+1} - S_j = P_j R_{j+1} - \mu^{j+1} - P_j + \mu^j.
			\]
			Conditioning on $R_1, \ldots, R_j$ and using the independence of $R_{j+1}$,
			\[
			\E(S_{j+1} - S_j \mid R_1, \ldots, R_j)
			= (\mu - 1)(P_j - \mu^j) = (\mu-1)S_j.
			\]
			Then,
			\[
			\E\bigl[(S_{j+1} - S_j)\,g(S_1, \ldots, S_j)\bigr]
			= (\mu - 1)\,\E\bigl[S_j\,g(S_1, \ldots, S_j)\bigr].
			\]
			Since $\E(S_j) = 0$, this reduces to
			$(\mu - 1)\Cov(S_j,\,g(S_1,\ldots,S_j))$.
			Now $R_1, \ldots, R_j$ are independent and hence associated. Each
			$S_k = \prod_{i=1}^k R_i - \mu^k$ is nondecreasing in $R_1, \ldots, R_k$ (and trivially
			in $R_{k+1}, \ldots, R_j$), so $g(S_1, \ldots, S_j)$ is also nondecreasing in
			$R_1, \ldots, R_j$. By association,
			$\Cov(S_j,\,g(S_1,\ldots,S_j)) \geq 0$, and since $\mu > 1$ the demimartingale
			inequality follows.
		\end{proof}
		
		\medskip
		
		\begin{remark}
			\label{rem:not-martingale-3}
			Since $\E(S_{n+1} \mid R_1, \ldots, R_n) = \mu S_n \neq S_n$ for $\mu \neq 1$, the
			process is not a martingale. The increment sequence need not be associated. To see this,
			take $R_1, R_2, \ldots$ i.i.d.\ with
			\[
			P(R = 0) = \tfrac{1}{5}, \quad P(R = 2) = \tfrac{3}{10}, \quad P(R = 4) = \tfrac{1}{2},
			\]
			so that $\mu = \frac{13}{5}$. The first two increments are
			\[
			X_1 = R_1(R_2 - 1) - \tfrac{104}{25}
			\quad\text{and}\quad
			X_2 = R_1 R_2(R_3-1) - \tfrac{1352}{125}.
			\]
			Define
			\[
			f(x) = I\!\left\{x > -\tfrac{54}{25}\right\},
			\qquad
			h(y) = I\!\left\{y > -\tfrac{2352}{125}\right\}.
			\]
			Direct calculation gives
			\[
			\E[f(X_1)] = \tfrac{11}{20},
			\qquad
			\E[h(X_2)] = \tfrac{89}{100},
			\qquad
			\E[f(X_1)h(X_2)] = \tfrac{11}{25},
			\]
			so
			\[
			\Cov(f(X_1),\,h(X_2))
			= \frac{11}{25} - \frac{11}{20} \cdot \frac{89}{100}
			= -\frac{99}{2000} < 0.
			\]
			Since $(X_1, X_2)$ is not associated, the increment sequence $\{X_n\}$ is not associated.
		\end{remark}
		
		\medskip
		
		\begin{remark}
			\label{rem:multiplicative-applications}
			The process $P_n = \prod_{i=1}^n R_i$ arises naturally in models of compounded investment
			returns, population growth in random environments, epidemic transmission across generations,
			and reliability degradation under successive random stress multipliers. In each of these
			settings, $S_n = P_n - \mu^n$ tracks the deviation of the multiplicative process from its
			expected trajectory. Two natural stopping times are
			\[
			\tau_b = \inf\{n \geq 1 : S_n > b\}
			\quad\text{and}\quad
			\tau_a = \inf\{n \geq 1 : S_n < -a\},
			\]
			corresponding respectively to the first time the process exceeds a target level and the
			first time a critical loss threshold is breached. The stopped process $S_\tau$ represents
			the state of the system at the moment the decision rule is activated, and the optional
			sampling results of this paper provide bounds on its expected value at that moment.
		\end{remark}

\medskip

		The preceding examples illustrate that stopping times arise in a natural manner in a variety of demimartingale models. We now turn to the main objective of the paper, namely the development of optional sampling inequalities for demi(sub)martingales under suitable monotonicity assumptions on the stopping rule. Optional sampling constitutes one of the fundamental principles of martingale theory and underlies many important results. However, despite the extensive literature on demimartingales (see \cite{PR2012} and references therein), corresponding optional sampling results appear to be unavailable. The aim of this paper is to address this gap by establishing expectation inequalities for stopped demi(sub)martingales under suitable monotonicity assumptions on the stopping rule.
		
		We begin in Section 2 by establishing comparison inequalities for demi(sub)martingales stopped at bounded stopping times. These finite-horizon results are subsequently extended in Section 3 to possibly unbounded stopping times under standard integrability conditions, leading to Doob-type optional sampling inequalities. Together, these results provide the first general optional sampling framework for demimartingales. Finally, Section 4 presents several applications of the developed theory. These include exponential and maximal inequalities for demi(sub)martingales, an almost sure bound for the maximal negative excursion of partial sums of associated random variables, and a collection of generalized Wald-type inequalities for positively associated sequences, including new inequalities for convex transformations of stopped sums.
	
	\medskip
	
	Throughout the paper, we denote $x\wedge y = \min\{x,y\}$ and $x^+ = \max\{0,x\}$ while $I\{A\}$ is used for the indicator function of the set $A$. Moreover, we use the convention that $\sum_{i=1}^{0} = 0$ where necessary while  the term \textit{demi(sub)martingale} means a demimartingale or a demisubmartingale sequence.
	
	First, we present here for convenience (and without proofs),  two theorems that are instrumental for the results obtained in this paper. The first can be found in \cite{C2003} and \cite{H2010} (check also Theorem 2.1.3 in \cite{PR2012}) while the second one is in \cite{PR2012} (see Theorem 2.1.4).
	
	\begin{theorem}
		\label{chr}Let the sequence $\left\{S_n, n \geq 1\right\}$ be a demisubmartingale or a demimartingale, $S_0=$ 0, and $\tau$ be a positive integer-valued random variable. Furthermore suppose that the indicator function $I\{\tau \leq j\}=h_j\left(S_1, \ldots, S_j\right)$ is a componentwise nonincreasing function of $S_1, \ldots, S_j$ for $j \geq 1$. Then the random sequence $\left\{S_j^*=S_{\tau \wedge j}, j \geq 1\right\}$ is a demisubmartingale.
	\end{theorem}
	
	\begin{theorem}
		\label{PR}Let the sequence $\left\{S_n, n \geq 1\right\}$ be a demisubmartingale and $\tau$ be a positive integer-valued random variable. Furthermore suppose that the indicator function $I\{\tau \leq j\}=h_j\left(S_1, \ldots, S_j\right)$ is a componentwise nonincreasing function of $S_1, \ldots, S_j$ for $j \geq 1$. Then, for any $1 \leq n \leq m$,
		
		$$
		\E\left(S_{\tau \wedge m}\right) \geq \E\left(S_{\tau \wedge n}\right) \geq \E\left(S_1\right).
		$$
		Suppose the sequence $\left\{S_n, n \geq 1\right\}$ is a demimartingale and the indicator function $I\{\tau \leq j\}=h_j\left(S_1, \ldots, S_j\right)$ is a componentwise nondecreasing function of $S_1, \ldots, S_j$ for $j \geq 1$. Then, for any $1 \leq n \leq m$,
		
		$$
		\E\left(S_{\tau \wedge m}\right) \leq \E\left(S_{\tau \wedge n}\right) \leq \E\left(S_1\right).
		$$
	\end{theorem}
	
	\section{Stopping time results for demi(sub)martingales}
	We begin with a basic inequality for demi(sub)martingales stopped at a bounded stopping time. 
	
	\begin{theorem}
		\label{theor12}Let $(S_n)_{n\geq 1}$ be a demi(sub)martingale and let $\tau$ be a bounded stopping time, i.e. there is $M\in\mathbb{N}$ such that $P(\tau \leq M) = 1$ such that $I\{\tau = k\}$ is a componentwise nondecreasing function of $S_1,\ldots, S_k$ for $k=1,\ldots, M-1$. Then,
		\begin{enumerate}
			\item[(a)] $\{S_\tau, S_M\}$ is a demisubmartingale
			\item [(b)]	 $\E(S_\tau)\leq \E(S_M)$.
		\end{enumerate}
	\end{theorem}	
	\begin{proof}
		Let $T_1 = S_{\tau}$ and $T_2 = S_M$ and let $f$ be a nonnegative nondecreasing function. Then,
		\begin{align*}
			&\E[(T_2-T_1)f(T_1)] = E[(S_M- S_\tau)f(S_\tau)] = \E[(S_M- S_\tau)f(S_\tau)I\{\tau\leq M\}]\\
			& = \sum_{k=1}^{M}  \E[(S_M- S_k)f(S_k)I\{\tau=k\}] \geq 0
		\end{align*}
		due to the demisubmartingale property. Thus, the pair $\{T_1, T_2\}$ forms a demisubmartingale sequence. The inequality in part (b) is due to the demisubmartingale property by choosing $f\equiv 1$.
	\end{proof}
	
	\medskip
	\noindent The corollary below follows as a direct consequence of Theorem \ref{theor12}.
	
	\begin{corollary}
		Let $(S_n)_{n\geq 1}$ be a demi(sub)martingale and let $\tau$ be a stopping time such that $I\{\tau = k\}$ is a componentwise nondecreasing function of $S_1,\ldots, S_k$. Then,
		\[
		\E(S_{\tau\wedge j})\leq \E(S_j), \quad j=1,2,\ldots.
		\]
	\end{corollary}
	\begin{proof}
		The result follows by applying the previous theorem since $\tau\wedge j$ is bounded by $j$.
	\end{proof}

	\begin{theorem}
	\label{th1}	Let $(S_n)_{n\ge1}$ be a demi(sub)martingale and let $\tau$ be a bounded stopping time such that
		\[
		\tau\le M \qquad \text{a.s.}
		\]
		Assume that, for every $i=2,\ldots,M$,
		\[
		I\{\tau<i\}
		\]
		is a componentwise nondecreasing function of $S_1,\ldots,S_{i-1}$. Then
		\[
		\E(S_\tau)\le \E(S_M).
		\]
	\end{theorem}
	\begin{proof}
		Since $\tau\le M$, we have
		\[
		S_M-S_\tau
		=
		\sum_{i=2}^{M}(S_i-S_{i-1})I\{\tau<i\}.
		\]
		Taking expectations gives
		\[
		\E(S_M-S_\tau)
		=
		\sum_{i=2}^{M}
		\E\left[(S_i-S_{i-1})I\{\tau<i\}\right].
		\]
		For each $i=2,\ldots,M$, the indicator $I\{\tau<i\}$ is nonnegative and, by assumption, componentwise nondecreasing in $S_1,\ldots,S_{i-1}$. Hence, by the demisubmartingale property,
		\[
		\E\left[(S_i-S_{i-1})I\{\tau<i\}\right]\ge0.
		\]
		Therefore,
		\[
		\E(S_M-S_\tau)\ge0,
		\]
		and consequently
		\[
		\E(S_\tau)\le \E(S_M).
		\]
	\end{proof}
	
	\begin{corollary}
	\label{cor1}	Let $(S_n)_{n\ge 1}$ be a demi(sub)martingale and let $\tau$ be a stopping time such that, for every $i=2,\ldots,j$,
		\[
		I\{\tau<i\}
		\]
		is a componentwise nondecreasing function of $S_1,\ldots,S_{i-1}$. Then
		\[
		\E(S_{\tau\wedge j})\le \E(S_j),
		\qquad j=1,2,\ldots.
		\]
	\end{corollary}
	\begin{proof}
		Apply Theorem \ref{th1} with $M=j$ and $\sigma=\tau\wedge j$. Then $\sigma\le j$ a.s. and for each
		$i=2,\ldots,j$,
		\[
		\{\sigma<i\}=\{\tau\wedge j<i\}=\{\tau<i\}.
		\]
		Indeed, since $i\le j$, if $\tau\wedge j<i$, then necessarily
		$\tau<i$; the converse is immediate. Hence
		\[
		I\{\sigma<i\}=I\{\tau<i\}, \qquad i=2,\ldots,j.
		\]
		By the assumption on $\tau$, it follows that $I\{\sigma<i\}$ is a
		componentwise nondecreasing function of $S_1,\ldots,S_{i-1}$ for every
		$i=2,\ldots,j$. Therefore, all assumptions of Theorem \ref{th1} are satisfied
		with $M=j$ and stopping time $\sigma$. Consequently,
		\[
		\E S_{\tau\wedge j}
		=
		\E S_\sigma
		\le
		\E S_j.
		\]
	\end{proof}

\begin{theorem}
\label{th2}	Let $(S_n)_{n\ge1}$ be a demi(sub)martingale. Let $r$ and $M$ be fixed integers such that
	\[
	1\le r<M,
	\]
	and let $\tau$ be a bounded stopping time satisfying
	\[
	r\le \tau\le M
	\qquad \text{a.s.}
	\]
	Assume that, for every $i=r+1,\ldots,M$,
	\[
	I\{\tau\ge i\}
	\]
	is a componentwise nondecreasing function of
	$S_1,\ldots,S_{i-1}$.
	Then $(S_r,S_\tau)$ is a demisubmartingale and therefore $\E(S_r)\leq \E(S_\tau).$
\end{theorem}
\begin{proof}
	Let $g$ be a nonnegative nondecreasing function. Since
	$r\le\tau\le M$, we have
	\[
	S_\tau-S_r
	=
	\sum_{i=r+1}^{M}(S_i-S_{i-1})I\{\tau\ge i\}.
	\]
	Therefore,
	\[
	(S_\tau-S_r)g(S_r)
	=
	\sum_{i=r+1}^{M}
	(S_i-S_{i-1})g(S_r)I\{\tau\ge i\}.
	\]
	For each $i=r+1,\ldots,M$, define
	\[
	F_i(S_1,\ldots,S_{i-1})
	=
	g(S_r)I\{\tau\ge i\}.
	\]
	Since $r<i$, the function
	\[
	(x_1,\ldots,x_{i-1})
	\mapsto g(x_r)
	\]
	is nonnegative and componentwise nondecreasing. By assumption,
	$I\{\tau\ge i\}$ is also componentwise nondecreasing in
	$S_1,\ldots,S_{i-1}$. Hence $F_i$ is a nonnegative
	componentwise nondecreasing function of
	$S_1,\ldots,S_{i-1}$. Applying the demisubmartingale property yields
	\[
	\E\Bigl[(S_i-S_{i-1})F_i(S_1,\ldots,S_{i-1})\Bigr]
	\ge0,
	\qquad i=r+1,\ldots,M.
	\]
	Summing over $i=r+1,\ldots,M$, we obtain
	\begin{align*}
		\E\bigl[(S_\tau-S_r)g(S_r)\bigr]
		&=
		\sum_{i=r+1}^{M}
		\E\Bigl[(S_i-S_{i-1})F_i(S_1,\ldots,S_{i-1})\Bigr]\geq0
	\end{align*}
	which confirms the demisubmartingale property for the pair $(S_r, S_\tau)$.
\end{proof}

\begin{remark}
	The monotonicity assumptions imposed above are not merely technical. They are satisfied by first-passage and threshold-type stopping times. For instance, let
	\[
	\tau_a=\inf\{n\geq 1:S_n\geq a\}\wedge M,
	\]
	where \(a\in\mathbb{R}\). Then, for every \(i=2,\ldots,M\),
	\[
	\{\tau_a<i\}
	=
	\left\{\max_{1\leq k\leq i-1}S_k\geq a\right\}.
	\]
	Hence
	\[
	I\{\tau_a<i\}
	=
	I\left\{\max_{1\leq k\leq i-1}S_k\geq a\right\}
	\]
	is a componentwise nondecreasing function of \(S_1,\ldots,S_{i-1}\). Therefore Theorem \ref{th1} applies and gives
	\[
	\E(S_{\tau_a})\leq \E(S_M).
	\]
	Similarly, Corollary \ref{cor1} gives
	\[
	\E(S_{\tau_a\wedge j})\leq \E(S_j),
	\qquad j\geq 1.
	\]
	The condition appearing in Theorem \ref{th2} is also satisfied by common choices of stopping times. For example, fix \(b\in\mathbb{R}\) and \(1\leq r<M\), and define
	\[
	\tau_b=\inf\{n\geq r:S_n<b\}\wedge M.
	\]
	Then, for every \(i=r+1,\ldots,M\),
	\[
	\{\tau_b\geq i\}
	=
	\{S_r\geq b,\ldots,S_{i-1}\geq b\}
	=
	\left\{\min_{r\leq k\leq i-1}S_k\geq b\right\}.
	\]
	Consequently,
	\[
	I\{\tau_b\geq i\}
	=
	I\left\{\min_{r\leq k\leq i-1}S_k\geq b\right\}
	\]
	is componentwise nondecreasing in \(S_1,\ldots,S_{i-1}\). Hence the assumptions of Theorem~2.5 are satisfied and the pair \((S_r,S_{\tau_b})\) is a demisubmartingale.
\end{remark}

	\section{Doob-type optional sampling theorems for demimartingales}

		The results of Section 2 compare a stopped process with a deterministic terminal value $S_M$ and apply to bounded stopping times. In many problems of interest, however, stopping times are unbounded, and a more reasonable comparison is with the initial value $S_1$ in the spirit of Doob's classical theorem. We establish such inequalities below under standard integrability conditions.

	\begin{theorem}
		\label{DoobOpt}	Let $(S_n)_{n\geq 1}$ be a demimartingale and let $\tau$  be a stopping time such that $I\{\tau\leq j\}$ is a componentwise nondecreasing function of $S_1,\ldots,S_j$. If any one of the following conditions 
		\begin{enumerate}
			\item [$(A_1)$] $\tau$ is bounded
			\item [$(A_2)$]  $(S_n)_{n\geq 1}$ is bounded and $\tau$ is finite a.s.
			\item [$(A_3)$] $\E(\tau)<\infty$ and $|S_k-S_{k-1}|\leq M$ for all $k\geq 2$ and for some $M>0$
		\end{enumerate}
		hold, then $\E(S_\tau) \leq \E(S_1)$.
	\end{theorem}
	\begin{proof}
		Notice first that if $\tau$ satisfies any of the conditions stated above, then $\tau <\infty$ a.s. Then, as $n\to \infty$ 
		\begin{equation}
			\label{conv1} \tau\wedge n \to \tau \quad \mbox{and}\quad
			S_{\tau\wedge n}\to S_\tau \quad \mbox{a.s.}. 
		\end{equation}
		Moreover, under the assumption that $I\{\tau\leq j\}$ is a componentwise nondecreasing function of $S_1,\ldots,S_j$ , Theorem \ref{PR} ensures that 
		\begin{equation}
			\label{ineq1}	\E(S_{\tau\wedge n}) \leq \E S_1.
		\end{equation}
		Let $(A_1)$ hold i.e., there is a constant $c$ such that $\tau\leq c$ a.s. Then, for all $ n>c$ we have that 
		$$\E(S_{\tau\wedge n}) = \E S_\tau\leq \E S_1$$
		where the last inequality follows by \eqref{ineq1}.
		
		\medskip
		
		\noindent Now, assume that $(A_2)$ is valid i.e. $\exists M<\infty$ such that $|S_n|<M$ for all $n$ and $\tau<\infty$ a.s. Based on this assumption we have that $|S_{\tau\wedge n}|<M$. Then,
		\begin{align*}
			&\E(S_{\tau}) = \E\left (\lim_{n\to\infty}S_{\tau\wedge n}\right ) = \lim_{n\to\infty}\E(S_{\tau\wedge n}) \leq  \lim_{n\to\infty}\E S_1 = \E S_1
		\end{align*}
		where the second equality follows by the dominated (bounded) convergence theorem while the inequality is again due to \eqref{ineq1}. 
		
		\medskip
		
		\noindent	Under assumption $(A_3)$ we start by noting that 
		\begin{align*}
			&|S_{\tau\wedge n} - S_1| = \left | \sum_{k=2}^{\tau\wedge n}(S_k-S_{k-1})\right |\leq  \sum_{k=2}^{\tau\wedge n}\left | (S_k-S_{k-1})\right | \leq (\tau\wedge n -1)M \leq (\tau -1)M.
		\end{align*}
		Notice that the RHS is independent of $n$ and since $\tau$ has finite expectation, $\E|S_{\tau\wedge n} - S_1| < \infty$. The desired result follows by employing once again the dominated (bounded) convergence theorem and \eqref{ineq1} i.e.
		\begin{align*}
			&\E(S_{\tau}-S_1) = \E\left  (\lim_{n\to\infty}(S_{\tau\wedge n}-S_1)\right ) = \lim_{n\to\infty}(\E(S_{\tau\wedge n}) - \E S_1)\leq 0.
		\end{align*}
		
		\medskip
		
	\end{proof}
	
	\noindent The next result strengthens Theorem \ref{DoobOpt} in the case where the demimartingale is nonnegative. 
	
	\begin{theorem}
		Consider a nonnegative demimartingale sequence $(S_n)_{n\geq 1}$  and let $\tau$  be a stopping time such that $\tau<\infty$ a.s. and  $I\{\tau\leq j\}$ is a componentwise nondecreasing function of $S_1,\ldots,S_j$. Then,
		\[
		\E S_\tau \leq \E S_1.
		\]
	\end{theorem}
	\begin{proof}
		The proof follows by applying Fatou's lemma. Particularly,
		\[
		\E S_\tau = \E\left (\lim_{n\to \infty} S_{\tau\wedge n}\right ) = \E\left (\liminf_{n\to \infty} S_{\tau\wedge n}\right ) \leq \liminf_{n\to \infty} \E  S_{\tau\wedge n} \leq \E S_1. 
		\]
	\end{proof}
	
	\noindent Theorem \ref{DoobOpt2} provides a counterpart to Theorem \ref{DoobOpt} for demisubmartingales when the stopping time satisfies a nonincreasing monotonicity condition. 
	
	\begin{theorem}
		\label{DoobOpt2}	Let $(S_n)_{n\geq 1}$ be a demisubmartingale and let $\tau$  be a stopping time such that $I\{\tau\leq j\}$ is a componentwise nonincreasing function of $S_1,\ldots,S_j$. If any one of the conditions $(A_1)$ - $(A_3)$ is true then, $\E(S_\tau) \geq \E(S_1)$.
	\end{theorem}
	\begin{proof}
		The results follow by applying similar steps as in Theorem \ref{DoobOpt} and taking into account that under the assumption on $I\{\tau\leq j\}$ we have that
		\[
		\E S_{\tau\wedge n} \geq \E S_1
		\]
		(by Theorem \ref{PR}).
	\end{proof}
	
	\section{Applications}
	\subsection{Maximal inequalities and almost sure bounds for associated sequences}
	
	The main result of this subsection is a lower-tail maximal inequality for demisubmartingales, obtained via a first-passage stopping time argument. This leads to an almost sure bound on the maximal negative excursion of partial sums of positively associated random variables. An exponential inequality is also derived as a further consequence of the optional sampling results of Section 3.
	
	\begin{theorem}
		\label{prob_min}	Let $(S_n)_{n\geq 1}$ be a demisubmartingale. Then, for any $\lambda >0$,
		\[
		P\left(\min_{1\leq k \leq n}S_k \leq -\lambda\right) \leq \dfrac{\E S_n^+-\E S_1}{\lambda}.
		\]
	\end{theorem}
	\begin{proof}
		For $\lambda>0$, define $\tau = \inf\{1\leq k \leq n: S_k\leq -\lambda\}\wedge n$. For this stopping time, $I\{\tau\leq j\}$ is a nonincreasing function of $S_1,\ldots, S_j$. Then, by Theorem \ref{DoobOpt2},
		$
		\E S_\tau \geq \E S_1.
		$
		Define $A = \{\tau = n, S_n \leq -\lambda\}$ (first hit occurs at $n$) and $B = \{\tau = n, S_n > -\lambda\}$ (no hit by time $n$). Then,   
		\begin{align*}
			&\E(S_\tau) = \E(S_\tau I\{\tau <n\})+\E(S_n I(A)) +\E(S_n I(B)). 
		\end{align*}
		On $\{\tau<n\}\cup A$ we have that $S_\tau \leq -\lambda$. Hence,
		\[
		\E(S_\tau I\{\tau <n\})\leq -\lambda P(\tau<n)
		\]
		while $\E(S_n I(A)) \leq -\lambda P(A)$. The latter two inequalities lead to
		\[
		\E(S_\tau I\{\tau <n\})+\E(S_n I(A)) \leq -\lambda P(\{\tau<n\}\cup A) = -\lambda P\left(\min_{1\leq k \leq n}S_k \leq -\lambda\right). 
		\]
		Moreover, $\E(S_nI(B)) \leq \E S_n^+$. Thus, the result follows by
		\[
		\E S_1 \leq \E S_\tau \leq -\lambda P\left(\min_{1\leq k \leq n}S_k \leq -\lambda\right) + \E S_n^+.
		\]
	\end{proof}
	
	\medskip
	
\noindent 	The next result is a consequence of Theorem  \ref{prob_min} and will lead to an almost sure convergence result. 
	
	\begin{corollary}
		\label{prob_as}Let $X_1, X_2, \ldots,$ be positively associated random variables with $\E X_k \geq0$ for all $k$. Then, for every $\lambda>0$
		\[
		P\left(\min_{1\leq k\leq n}S_k \leq -\lambda \right)\leq \dfrac{\E S_n^+}{\lambda}\leq \dfrac{\sum_{i=1}^{n}\E X_i^+}{\lambda}.
		\]
	\end{corollary}

	\begin{corollary}
	\label{corconv}	Let $X_1, X_2, \ldots,$ be positively associated random variables with $\E X_k \geq0$ for all $k$ with $\sup_i \E X_i ^+ <\infty$. Then, for every $\varepsilon>0$,
		\[
		\dfrac{\displaystyle\max_{1\leq k\leq n}(-S_k)^+}{n(\log n)^{1+\varepsilon}}\to 0, \quad a.s., \quad n\to\infty. 
		\]
	\end{corollary}
	\begin{proof}
		Let $\sup_i \E X_i ^+ \equiv M<\infty$ and notice that $D_n =\displaystyle \max_{1\leq k\leq n}(-S_k)^+ = (-\displaystyle\min_{1\leq k \leq n}S_k)^+$. 
		The latter equality is obtained as follows:
		Set $m = \min_{1\leq k \leq n} S_k$ and pick $k^*$ such that $S_{k^*} = m$. Then,
		\[
		\max_{1\leq k\leq n}(-S_k)^+ = \max\{(-S_1)^+, (-S_2)^+,\ldots, (-S_n)^+\} \geq (-S_{k^*} )^+ = (-\min_{1\leq k \leq n}S_k)^+.
		\]				
		On the other hand, for any $k= 1,2, \ldots, n$ we have that $S_k \geq m$ which leads to $(-S_k)^+ \leq (-m)^+$ since $h(x) = x^+$ is a nondecreasing function. Thus, $\max_{1\leq k\leq n}(-S_k)^+ \leq (-m)^+$.
		
		\medskip
		
		\noindent Moreover, notice that for any $\lambda>0$ we have that 
		$$P \left (\left (-\displaystyle\min_{1\leq k \leq n}S_k\right )^+ \geq \lambda\right ) =P\left (-\displaystyle\min_{1\leq k \leq n}S_k \geq \lambda\right ).  $$
		
		\medskip
		
		\noindent From Corollary \ref{prob_as}, for every $\lambda>0$ we have that 
		\[
		P(D_n \geq \lambda) = P\left (-\displaystyle\min_{1\leq k \leq n}S_k \geq \lambda\right ) \leq  \dfrac{\sum_{i=1}^{n}\E X_i^+}{\lambda} \leq \dfrac{nM}{\lambda}.
		\]
		For $ m\in \mathbb{N}$, let $n_m = 2^m$, and fix an arbitrary $c>0$. Set $\lambda_m = c2^m m^{1+\varepsilon}$.  Then,
		\[
		P(D_{2^m}\geq \lambda_m) \leq \dfrac{M}{c m^{1+\varepsilon}}
		\]
		which ensures that $\sum_{m=1}^{\infty}P(D_{2^m}\geq \lambda_m)<\infty$ and by the Borel-Cantelli lemma we have that
		\[
		P((D_{2^m}\geq c2^m m^{1+\varepsilon}\quad\mbox{i.o}) = 0.
		\]
		Therefore, for every fixed $c>0$,
		\[
		\limsup_{m\to\infty}\dfrac{D_{2^m} }{2^m m^{1+\varepsilon}}\leq c\quad {\rm a.s.}
		\]
		Since $c>0$ is arbitrary, taking $c = 1/r,\, \, r=1,2, \ldots$ and intersecting the corresponding probability-one events, we obtain
		\[
		\limsup_{m\to\infty}\dfrac{D_{2^m} }{2^m m^{1+\varepsilon}}=0 \quad {\rm a.s.}
		\]
		which allows for the conclusion $$\dfrac{D_{2^m} }{2^m m^{1+\varepsilon}}\to 0 \quad {\rm a.s.}$$ 
		Now let \(n\in\mathbb N\) and choose \(m\) such that
\[
2^m\le n<2^{m+1}.
\]
Since \(D_n\) is nondecreasing in \(n\), we have
\[
D_n\le D_{2^{m+1}}.
\]
Therefore
\[
0\le
\frac{D_n}{n(\log n)^{1+\varepsilon}}
\le
\frac{D_{2^{m+1}}}{2^{m+1}(m+1)^{1+\varepsilon}}
\cdot
\frac{2^{m+1}(m+1)^{1+\varepsilon}}
{n(\log n)^{1+\varepsilon}}.
\]
The first factor tends to \(0\) almost surely by the dyadic convergence just proved. Moreover,
since \(2^m\le n<2^{m+1}\),
\[
\frac{2^{m+1}(m+1)^{1+\varepsilon}}
{n(\log n)^{1+\varepsilon}}
\le
2\left(\frac{m+1}{m\log 2}\right)^{1+\varepsilon},
\]
which is bounded for all sufficiently large \(m\). Hence
\[
\frac{D_n}{n(\log n)^{1+\varepsilon}}\to0
\qquad \text{a.s.}
\]
\end{proof}

\begin{remark}
	The quantity $D_n$	represents the maximal negative excursion of the partial-sum process up to time \(n\) and it arises in sequential analysis, where it measures the largest cumulative negative deviation of a monitored process, and in ruin theory, where it represents the maximal deficit of an insurer up to time n. Corollary \ref{corconv} provides an almost sure asymptotic bound on this quantity under positive dependence.
\end{remark}


	\subsection{Generalized Wald-type inequalities for associated random variables} 
	Wald-type identities and inequalities play a fundamental role in probability theory as they relate characteristics of a randomly stopped sum to those of the underlying stopping time. The classical Wald identity for i.i.d. random variables was established by Wald (1944) in \cite{W1944} and has since become a cornerstone of sequential analysis and renewal theory. Extensions to independent but not necessarily identically distributed random variables were obtained by Robbins and Samuel (1966), while Chow, Robbins and Teicher (1965) developed higher-order moment identities and inequalities for randomly stopped sums (see \cite{RS1966} and \cite{CRH1965} respectively). More recently, de la Peña and Zamfirescu (2002) in \cite{PZ2002} established Wald-type identities for martingale difference sequences under suitable domination assumptions. A comprehensive study of Wald's identities and inequalities and their applications can be found in \cite{B2013}. Despite this extensive literature, corresponding results for positively associated random variables appear to be largely unavailable. Theorems \ref{WaldLinear} and \ref{GenWald} address this gap by establishing Wald-type inequalities for positively associated random variables and, more generally, for their nonlinear convex transformations.
	
	\medskip
	\noindent For all the results that follow, $S_n = \sum_{i=1}^{n}X_i$. Moreover, for a bounded stopping time $\tau$ and for a deterministic sequence
	$(a_i)_{i\ge1}$, we define
	\[
	\sum_{i=1}^{\tau} a_i
	:=
	\sum_{i=1}^{\infty} a_i I\{\tau\ge i\}.
	\]
	Since $\tau$ is bounded, say $\tau\le M$ a.s., this is a finite sum
	almost surely and
	\[
	\sum_{i=1}^{\tau} a_i
	=
	\sum_{i=1}^{M} a_i I\{\tau\ge i\}.
	\]
	Consequently,
	\[
	E\!\left(\sum_{i=1}^{\tau} a_i\right)
	=
	\sum_{i=1}^{M} a_i P(\tau\ge i).
	\]
	\begin{theorem}
		\label{WaldLinear}Let $X_1,X_2,\ldots$ be integrable positively associated random variables and let $\tau$ be a bounded stopping time.
		\begin{enumerate}
			\item [(a)]	If $I\{\tau\le j\}$ is a componentwise nondecreasing function of
			$S_1,\ldots,S_j$ for every $j\ge1$, then
			\[
			\E(S_\tau)
			\le
			\E\left(\sum_{i=1}^{\tau}\E(X_i)\right).
			\]
			\item [(b)]	If $I\{\tau\le j\}$ is a componentwise nonincreasing function of
			$S_1,\ldots,S_j$ for every $j\ge1$, then
			\[
		\E(S_\tau)
			\ge
			\E\left(\sum_{i=1}^{\tau}\E(X_i)\right).
			\]
		\end{enumerate}
	\end{theorem}
	\begin{proof}
		Define
		\[
		H_n=S_n-\sum_{i=1}^n\E(X_i)
		=
		\sum_{i=1}^n\{X_i-\E(X_i)\}.
		\]
		Notice that $X_i-\E(X_i)$ are mean-zero positively associated random variables. Therefore $(H_n)$ is a demimartingale and a demisubmartingale with $\E(H_1) = 0$. Moreover, $S_i = H_i + \sum_{k=1}^{i}\E(X_i)$ so the monotonicity of the indicator function transfers from $S_1, \ldots, S_j$ to $H_1, \ldots, H_j$. Thus, the first part of the theorem is obtained by applying Theorem \ref{DoobOpt} while the second part derives from Theorem \ref{DoobOpt2}.
		\end{proof}
		
		\begin{theorem}
		\label{GenWald}	Consider $(X_n)_{n\ge1}$ a sequence of nonnegative positively associated
			random variables. Let $g:[0,\infty)\to\mathbb{R}$ be convex and
			strictly increasing with $g(0)=0$, such that
			\[
			\E|g(X_n)|<\infty,\qquad n\ge1.
			\]
			Let $\tau$ be a bounded stopping time and assume that 
			\[
			I\{\tau\le j\}
			\]
			is a componentwise nonincreasing function of
			$S_1,\ldots,S_j$ for every $j\ge1$. Then,
			\[
			\E(g(S_\tau))
			\ge
			\E\left(\sum_{i=1}^{\tau}\E(g(X_i))\right).
			\]
			\end{theorem}
	\begin{proof}
		Let $j\ge1$ and let $f$ be a nonnegative componentwise nondecreasing
		function. Since $g$ is convex and $g(0)=0$, for all $x,y\ge0$, it follows that $g$ is superadditive on $[0,\infty)$, that is,
		\[
		g(x+y)-g(x)\ge g(y).
		\]
		Therefore,
		\[
		g(S_{j+1})-g(S_j)
		=
		g(S_j+X_{j+1})-g(S_j)
		\ge
		g(X_{j+1}).
		\]
		Define $H_n = g(S_n) - \sum_{i=1}^n \E(g(X_i))$; then,
		\[
		H_{j+1}-H_j
		=
		g(S_{j+1})-g(S_j)-\E(g(X_{j+1}))
		\ge
		g(X_{j+1})-\E(g(X_{j+1})).
		\]
		Thus,
		\[
		\E[(H_{j+1}-H_j)f(H_1,\ldots,H_j)]
		\ge
		\E[(g(X_{j+1})-\E(g(X_{j+1})))f(H_1,\ldots,H_j)].
		\]
		Since $g$ is nondecreasing, each $H_i$ is an increasing function of
		$X_1,\ldots,X_i$. Hence $f(H_1,\ldots,H_j)$ is a nondecreasing function
		of $X_1,\ldots,X_j$. Moreover, $g(X_{j+1})$ is a nondecreasing function
		of $X_{j+1}$. By positive association,
		\[
		\E[g(X_{j+1})f(H_1,\ldots,H_j)]
		\ge
		\E(g(X_{j+1}))\E[f(H_1,\ldots,H_j)].
		\]
		Therefore,
		\[
		\E[(g(X_{j+1})-\E(g(X_{j+1})))f(H_1,\ldots,H_j)]\ge0.
		\]
		Consequently,
		\[
		\E[(H_{j+1}-H_j)f(H_1,\ldots,H_j)]\ge0
		\]
		which proves the  demisubmartingale property for $(H_n)_{n\ge1}$. Since $g$ is strictly increasing, it admits an increasing inverse on its
		range. Moreover,
		\[
		S_i
		=
		g^{-1}\left(H_i+\sum_{k=1}^i \E(g(X_k))\right).
		\]
		Thus $S_i$ is a nondecreasing function of $H_i$. Therefore, if
		$I\{\tau\le j\}$ is componentwise nonincreasing in
		$S_1,\ldots,S_j$, then it is also componentwise nonincreasing in
		$H_1,\ldots,H_j$. Since $\tau$ is bounded, Theorem \ref{DoobOpt2} applies and yields
		\[
		\E(H_\tau)\ge \E(H_1) = 0
		\]
	which is written as
	\[
	\E\left(g(S_\tau)-\sum_{i=1}^{\tau}\E(g(X_i))\right)\ge0.
	\]
	\end{proof}
	
	\begin{corollary}
		Let $(X_n)_{n\ge1}$ be a sequence of nonnegative identically distributed positively associated
		random variables. Consider $g:[0,\infty)\to\mathbb{R}$ be convex and
		strictly increasing with $g(0)=0$, such that
		\[
		\E|g(X_n)|<\infty,\qquad n\ge1.
		\]
		Let $\tau$ be a bounded stopping time and assume that 
		\[
		I\{\tau\le j\}
		\]
		is a componentwise nonincreasing function of
		$S_1,\ldots,S_j$ for every $j\ge1$. Then,
		\[
		\E(g(S_\tau))
		\ge
		\E(g(X_1))\E\tau.
		\]
	\end{corollary}

		\begin{remark}
			Different choices of the convex function $g$ in Theorem \ref{GenWald}  lead to different forms of nonlinear stopped inequalities. For example:
			
			\begin{itemize}
				\item[(a)] Choosing
				\[
				g(x)=x^p,\qquad p\ge1,
				\]
				yields
				\[
				\E(S_\tau^p)\ge \E\left(\sum_{i=1}^{\tau}\E(X_i^p)\right).
				\]
				In the identically distributed case this becomes
				\[
				\E(S_\tau^p)\ge \E(X_1^p)\E(\tau).
				\]
				
				\item[(b)] Choosing
				\[
				g(x)=e^{\theta x}-1,\qquad \theta>0,
				\]
				leads to the exponential-type inequality
				\[
				\E(e^{\theta S_\tau}-1)
				\ge
				\E\left(
				\sum_{i=1}^{\tau}
				\bigl(\E(e^{\theta X_i})-1\bigr)
				\right).
				\]
				
				\item[(c)] Choosing
				\[
				g(x)=x\log(1+x),\qquad x\ge0,
				\]
				gives
				\[
				\E(S_\tau\log(1+S_\tau))
				\ge
				\E\left(
				\sum_{i=1}^{\tau}
				\E(X_i\log(1+X_i))
				\right).
				\]
			\end{itemize}
		
		\end{remark}

		\subsection*{Acknowledgement}
	Work of the second author is supported by the scheme "INSA Honorary Scientist" at the CR Rao Advanced Institute of Mathematics, Statistics  and Computer Science, Hyderabad, India.

\end{document}